# The Seventeen Secondary Elements of Pythagorean Triangles

Konstantine Zelator

1. Introduction

Almost every introductory text in number theory has a short chapter on Diophantine equations with a short section on Pythagorean triangles or triples; these are right triangles with integer sides or edge lengths. Even the reader who knows very little about the known properties of Pythagorean triples must have come across the famous triple (3,4,5) and some of its non-primitive multiples $(3\delta, 4\delta, 5\delta)$; where $\delta$ is a positive integer (and very likely the next primitive triple (5,12,13) as well). But let us state the usual formulas that generate and describe all Pythagorean triangles. If $\alpha, \beta, \gamma$ are the integer side lengths of a Pythagorean triangle, with $\alpha$ being the hypotenuse, then without loss of generality,

$$\boxed{\beta = \delta \cdot (2mn),\ \gamma = \delta \cdot (m^2 - n^2),\ \alpha = \delta \cdot (m^2 + n^2) \text{ for some positive integers } \delta, m, n \\ \text{such that } (m, n) = 1 \text{ and } m + n \equiv 1 \pmod{2}, m > n} \quad (1)$$

The above parametric formulas can easily be found in number theory books. For instance, see reference [1] and [6].

As the reader can readily verify, if $\alpha, \beta, \gamma$ are given by formulas (1), then an easy computation establishes that $\beta^2 + \gamma^2 = \alpha^2$. For a student taking a first course in number theory, more work is needed to establish the converse (some background in elementary number theory is needed for that): Namely, that if the three positive integers $\alpha, \beta, \gamma$ satisfy $\beta^2 + \gamma^2 = \alpha^2$, then they must be represented by formulas (1).

Let us remark here, that it is standard notation in number in number theory for the symbol $(m, n)$ to represent the greatest common divisor between $m$ and $n$; if $(m, n) = 1$, we say that the integers $m$ and $n$ are relatively prime. Also, if $\delta = 1$ in (1), the Pythagorean triangle is called primitive. Also note that since for every integer $x$ either $x \equiv 0 \pmod{2}$ (i.e, $x$ is even) or $x \equiv 1 \pmod{2}$ (i.e, $x$ is odd), then condition $m + n \equiv 1 \pmod{2}$, necessitates that one of $m, n$ is even, the other odd (they have different parities).

Some readers must undoubtedly be very familiar not only with the derivation of formulas (1), but also with a number of properties of Pythagorean triangles. W. Sierpinski's book "Elementary Theory of Numbers", first published in 1964 $(see [1])$, is a treasure trove of material on Pythagorean triples. In addition, the book " Recreations in the Theory of Numbers", by Albert H.Bieler (see [2]), also has a wealth of material on the subject. On the ever fascinating subject of Pythagorean triangles there are many solved problems for the student to learn from; and at least as many, a plethora one might say, of unsolved problems to this day; not counting new ones that have been posed since. Nevertheless, out of the various textbooks in number theory this author has studied or glanced through, only one has a bit of material on the secondary elements of a Pythagorean triangle, namely a couple of exercises involving the internal radius $\rho$ (that is, the radius of the inscribed circle). This is the book " Number Theory with Applications" , by James A.Anderson & James M.Bell (see [3]). In L.E. Dickson's voluminous and monumental History of the Theory of Numbers, 26 printed pages are devoted to Pythagorean triangles and right triangles with rational side lengths in general (see reference [5]).However, there is less than a page of information on the secondary elements of Pythagorean triangles. (These include medians, the internal and external angle bisectors.)

Specifically on page 188 of Dickson's book, the Green mathematician Diaphantus (who lived in Alexandria, Egypt, circa 150-250 AD, the entire area of number theory, Diophantine equations, was named after him, to honor him as the first historically known mathematician in the last 2000 years, to under take a more systematic approach to solving equations in integers or rational numbers) is reported as having found the Pythagorean tripe (28,96,100), which has an angle bisector with exact length 35. (See last section (Numerical Examples ) of this paper-Family 4).Also, on ;page 189 (in Dickson's book), J.Kersey is mentioned as having discovered a family of right triangles with rational side lengths. Specifically, he showed that if $b, p, h$ (hypotenuse) are the rational sides of a right triangle, then the right triangle with side lengths, $|\overline{AC}| = p(p^2 + b^2), |\overline{AB}| = p(p^2 - b^2), |\overline{BC}| = p(2bp)$, has an angle bisector $AD$ (of the angle A) with length $|\overline{AD}| = h(p^2 - b^2)$. Finally, E. Thrriere' is reported as having discovered a specific rational right triangle with rational interior and exterior bisectors of one of the angles.

The purpose of this paper is to pursue the following questions: which of the secondary elements are integers? Which are just rational numbers? Let us list the secondary elements of a triangle (not just Pythagorean) with side lengths $\alpha, \beta, \gamma$:

1) R = the radius of the circumscribed circle.
   $\rho$ : the radius of the inscribed circle.
2) $h_\alpha, h_\beta, h_\gamma$ : the lengths of the three heights.
3) $\delta_\alpha, \delta_\beta, \delta_\gamma$ : the lengths of the three internal bisectors; that is, the lengths of the bisectors of the triangle's internal angles.
4) $d_\alpha, d_\beta, d_\gamma$ : the lengths of the three bisectors of the triangles' external angles.
5) $\rho_\alpha, \rho_\beta, \rho_\gamma$ : these are the radii of the three circles; each being tangential to the tree (straight) lines containing the edges of the triangle, but lying outside the triangle.
6) $\mu_\alpha, \mu_\beta, \mu_\gamma$ : the lengths of the triangle's three medians.

Below we list the findings of this paper. A few of them are entirely obvious, by inspection. Some other ones may not be immediately obvious, but they easily emerge after a straightforward calculation. It is a bit more difficult to establish the precise conditions for the rationality of $\delta_\beta$ and $d_\beta$; and it is challenging to prove the irrationality of $\mu_\gamma$ and $\mu_\beta$, in a complicate proof that involves the method of infinite descent.

Here are the findings:

a) The internal radius $\rho$ is an integer, in any Pythagorean triangle. This is the result of a straight forward computation. On the other hand, since $R = \dfrac{\alpha}{2}$, it is obvious that R is either an integer, or half an odd integer.

b) Since $h_\beta = \gamma, h_\gamma = \beta$, it is obvious that the two heights $h_\beta$ and $h_\gamma$ are always integers. On the other hand, $h_\alpha$ is always rational and it is integral only in certain kinds of non-primitive Pythagorean triangles.

c) The internal bisectors $\delta_\alpha$ and $\delta_\gamma$ are always irrational; the internal bisector $\delta_\beta$ can be integer only in certain kinds of non-primitive triangles; it is sometimes rational, sometimes irrational; precise conditions are given.

d) The external bisectors $d_\alpha$ and $d_\gamma$ are always irrational; $d_\beta$ can be rational, integral or irrational; precise conditions are given.

e) The tree radii $\rho_\alpha, \rho_\beta, \rho_\gamma$ are always integers.

f) Obviously $\mu_\alpha = R = \dfrac{\alpha}{2}$; so, $\mu_\alpha$ is either an integer or half an odd integer. The other two, $\mu_\beta$ and $\mu_\gamma$ are always irrational, but it is hard to establish their irrationality.

In fact, all Pythagorean triangles in ( c ) with $\delta_\beta$ rational can be parametrically described in three families; and we offer parametric formulas in two families of such triangles, describing all Pythagorean triangles in ( d ) with $d_\beta$ being rational.

2. **Illustrations**

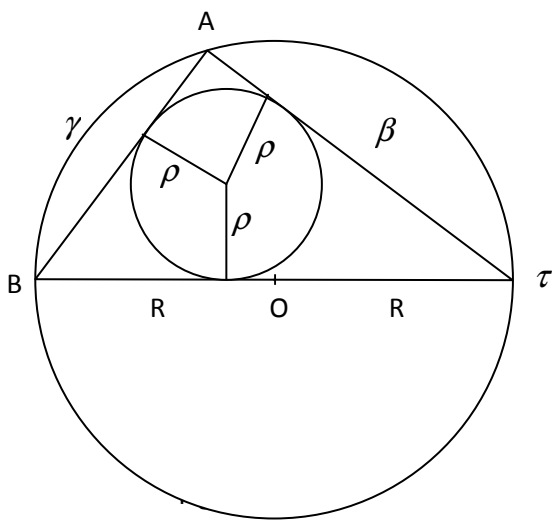

Figure 1

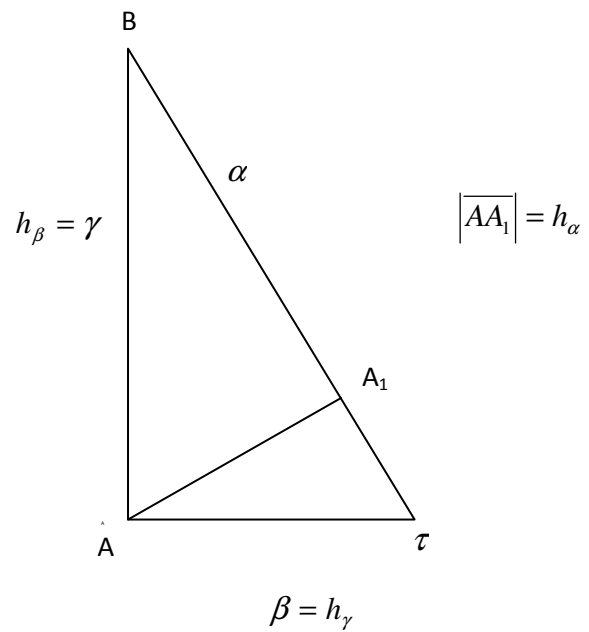

$\beta = h_\gamma$

Figure 2

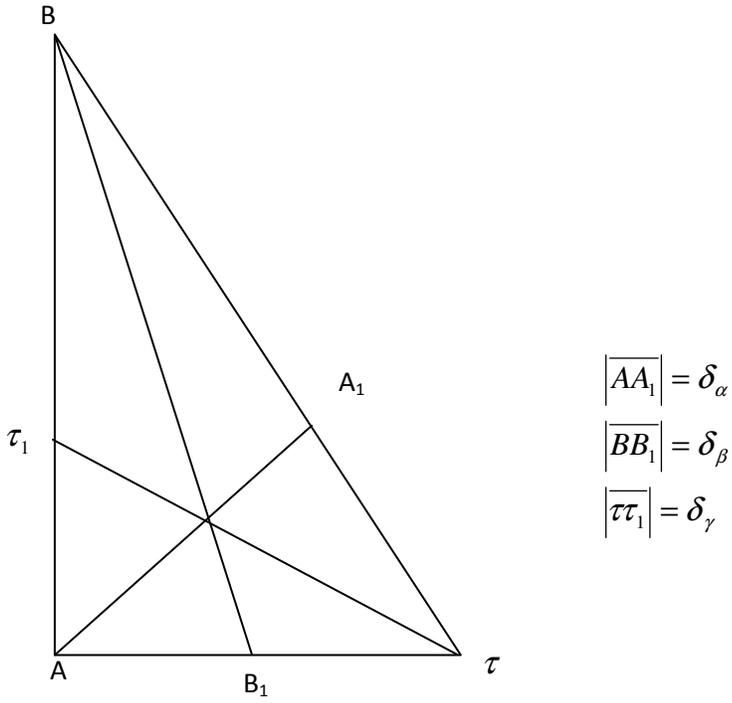

$$|\overline{AA_1}| = \delta_\alpha$$
$$|\overline{BB_1}| = \delta_\beta$$
$$|\overline{\tau\tau_1}| = \delta_\gamma$$

Figure 3

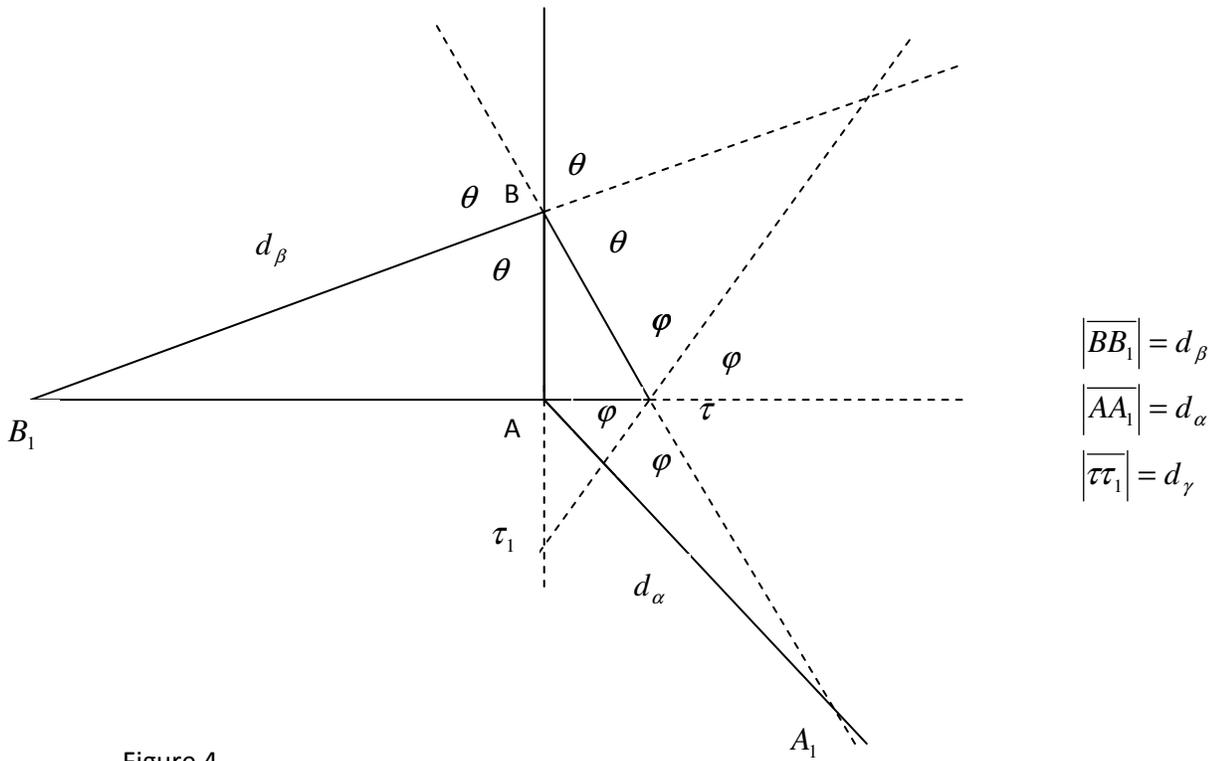

$$|\overline{BB_1}| = d_\beta$$
$$|\overline{AA_1}| = d_\alpha$$
$$|\overline{\tau\tau_1}| = d_\gamma$$

Figure 4

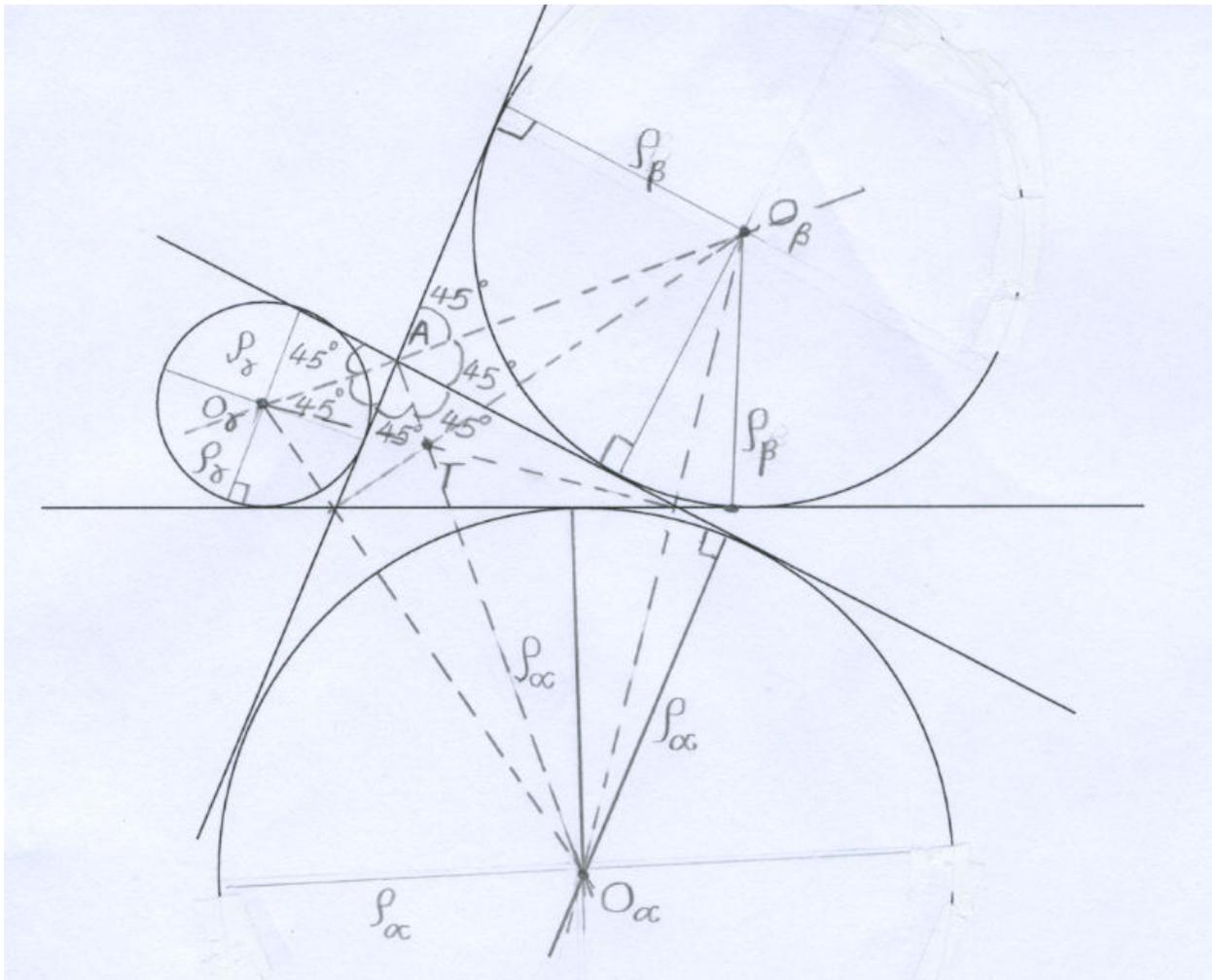

Figure 5

T is the point of intersection of the lines $O_\alpha A, O_\beta B, O_\gamma T$ of the angles $\hat{A}, \hat{B}, \hat{T}$ respectively. Also, these three line segments are the heights in the triangle $O_\alpha \overset{\Delta}{O}_\beta O_\gamma$. The points $O_\alpha, O_\beta, O_\gamma$ are the centers of the three exterior tangential circles.

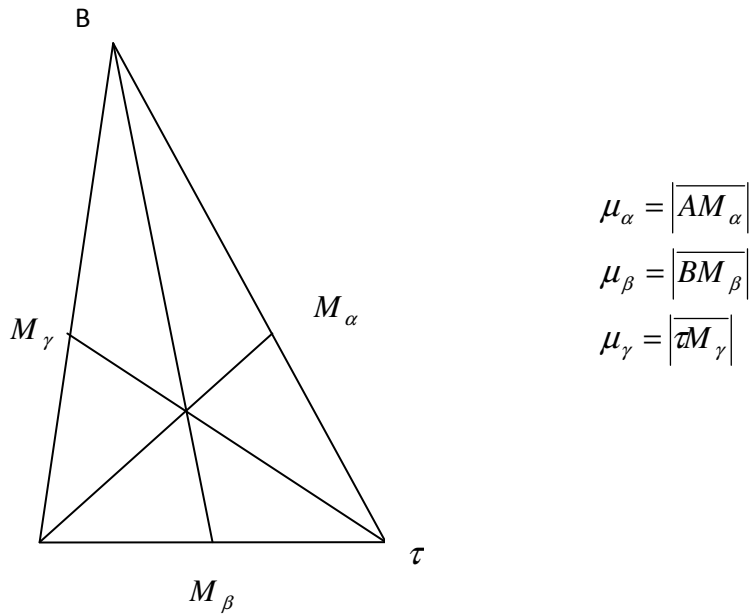

Figure 6

# 3 Triangle Formulas and a Few Facts from Number Theory

Some of the formulas listed below are well known in the literature of geometry and trigonometry, but not widely known among a general mathematical audience. However, a student typically makes little or no use of most of these formulas, in his/her high school and early college or university years.

For more information, refer to references [7] and [8]. For a triangle of side lengths $\alpha, \beta, \gamma$, these formulas are as follows:

Semi-perimeter $s = \dfrac{\alpha + \beta + \gamma}{2}$

Area $A = \sqrt{s(s-\alpha)(s-\beta)(s-\gamma)}$ (also known as Herson's formula)

For a right triangle with hypotenuse length $\alpha$, this becomes simplified to $A = \dfrac{\beta\gamma}{2}$

$(F_2)$: Heights: $h_\alpha = \dfrac{2A}{\alpha}, h_\beta = \dfrac{2A}{\beta}, h_\gamma = \dfrac{2A}{\gamma}$

When $\alpha$ is a hypotenuse length, these formulas become simplified to $h_\alpha = \dfrac{\beta\gamma}{\alpha}, h_\beta = \gamma, h_\gamma = \beta$

$(F_1)$: Radii of the circumscribed and inscribed circles: $R = \dfrac{\alpha\beta\gamma}{4A}, \rho = \dfrac{A}{s}$

When the triangle is a right one with $\alpha$ the hypotenuse length, $R = \dfrac{\alpha}{2}, \rho = \dfrac{\beta\gamma}{\alpha + \beta + \gamma}$

$(F_3)$: Internal bisectors $\delta_\alpha, \delta_\beta, \delta_\gamma$

$$\delta_\alpha = \dfrac{2\beta\gamma}{\beta+\gamma} \cdot \sqrt{\dfrac{s(s-a)}{\beta\gamma}}, \delta_\beta = \dfrac{2\alpha\gamma}{\alpha+\gamma} \cdot \sqrt{\dfrac{s(s-\beta)}{\alpha\gamma}}, \delta_r = \dfrac{2\alpha\beta}{\alpha+\beta} \cdot \sqrt{\dfrac{s(s-\gamma)}{\alpha\beta}}$$

$(F_4)$: External bisectors $d_\alpha, d_\beta, d_\gamma$

$$d_\alpha = \dfrac{2\beta\gamma}{|\beta-\gamma|} \cdot \sqrt{\dfrac{(s-\beta)(s-\gamma)}{\beta\gamma}}, d_\beta = \dfrac{2\alpha\gamma}{|\alpha-\gamma|} \cdot \sqrt{\dfrac{(s-\alpha)(s-\gamma)}{\alpha\gamma}},$$

$$d_\gamma = \dfrac{2\alpha\beta}{|\alpha-\beta|} \cdot \sqrt{\dfrac{(s-\alpha)(s-\beta)}{\alpha\beta}}$$

Note here that when a triangle is isosceles but not equilateral, exactly one of the external bisectors will be undefined; it will have infinite length so to speak. In such a case the (straight) line that contains that external bisector, will be parallel to the side opposite that bisector. When a triangle is equilateral, of course, all three $d_\alpha, d_\beta, d_\gamma$, will be undefined. Also, even though a right triangle can be isosceles, no Pythagorean triangle can be isosceles.

$(F_5)$: The radii of the three external tangential circles:

$$\rho_\alpha = \frac{A}{s-\alpha}, \rho_\beta = \frac{A}{s-\beta}, \rho_\gamma = \frac{A}{s-\gamma}$$

$(F_6)$: The three medians:

$$\mu_\alpha = \frac{\sqrt{2(\beta^2+\gamma^2)-\alpha^2}}{2}, \mu_\beta = \frac{\sqrt{2(\alpha^2+\gamma^2)-\beta^2}}{2}, \mu_\gamma = \frac{\sqrt{2(\alpha^2+\beta^2)-\gamma^2}}{2}$$

In the case for a right triangle with hypotenuse length $\alpha$, the first formula simplifies to $\mu_\alpha = \frac{\alpha}{2}$.

We also list seven basic facts from number theory.

### Seven Basic Facts form number theory

**Fact 1:** Every nonempty subset of the set of positive integers or natural numbers, must contain a least element.

**Fact 2:** For every prime $\rho \equiv 3 \pmod 4$, and any two integers a and b, the sum of squares $a^2+b^2$ is divisible by $\rho$ if, and only if, both a and b are divisible by $\rho$. In other words for every prim $\rho \equiv 3 \pmod 4$, $a^2+b^2 \equiv 0 \pmod p \Leftrightarrow a \equiv b \equiv \pmod p$

**Fact 3:** Let $a,b,c$ be integers such that $a$ divides the product $bc$; then, if $(a,b)=1$, $a$ must be a divisor of $c$.

**Fact 4:** If $a,b,c$ are integers such that $(a,b)=1=(a,c)$, then $(a,bc)=1$. More generally, if each of the integers $a_1, a_2,...,a_k$ is relatively prime to each of the integers $b_1, b_2,...,b_l$, then the product $a_1, a_2,...,a_k$ is relatively prime to the product $b_1, b_2,...,b_l$. As a result, if $a,b$ are integers with $(a,b)=1$; then $(a^k,b^l)=1$, for any positive integer exponents $k$ and $l$.

**Fact 5:** If two positive integers $a$ and $b$ divide each other, then $a=b$.

**Fact 6:** If (the integers $a$ and $b$ satisfy) $(a,b)=1$, then $(a+b,c-b)=1$ or 2 depending on whether $a$ and $b$ have different parity $(i.e, a+b \equiv 1 \pmod 2)$ or they are both odd.

**Fact 7:** If (the integers $a$ and $b$ satisfy) $(a,b)=1$, then $(a \pm b, ab)=1$. Also $(a \pm b, 2ab)=1$, if $a$ and $b$ have different parity; while $(a \pm b, 2ab)=2$, if both $a$ and $b$ are odd. More generally, $(a^{\varepsilon_1} \pm b^{\varepsilon_2}, a^{\varepsilon_1} \cdot b^{\varepsilon_2})=1$, for any positive integers $\varepsilon_1, \varepsilon_2$

We remark her, that **Fact 1** is exiomatic in nature; the first part of **Fact 2** as well as **Facts 3** through **7**, would typically be assigned as exercises in the first half of a first course in elementary number theory.

## 4   Computations and Conclusions

On $R, \rho, h_\alpha, h_\beta, h_\gamma, \delta_\alpha, \delta_\beta, \delta_\gamma, d_\alpha, d_\beta, d_\gamma, \rho_\alpha, \rho_\beta, \rho_\gamma$ and $\mu_\alpha$ ( For $\mu_\beta, \mu_\gamma$, see next section).

We now apply formulas $(F_1)$ through $(F_5)$ and the first formula in $(F_6)$ in the case of a pythagorean triangle. Using these formulas in conjunction with formulas (1) we arrive at the following simplified results.

(i) $R = \dfrac{\delta(m^2 + n^2)}{2}, \rho = \delta \cdot n \cdot (m - n)$

Since $m^2 + n^2$ is an odd integer, R will be an integer exactly when $\delta$ is an even number.

### Conclusion:

1) $\rho$ is always an integer.

2) R is always a rational number; it is an integer only in non-primitive pythagorean triangles with $\delta \equiv 0 (\mod 2)$; otherwise it is half an odd integer.

(ii) $h_\alpha = \dfrac{\delta \cdot (2mn)(m^2 - n^2)}{m^2 + n^2}, h_\beta = \delta \cdot (m^2 - n^2), h_\gamma = \delta \cdot (2mn)$

Since $(m, n) = 1$ and $m, n$ have different parity, **Facts 6 and 7** tell us that
$(m^2 + n^2, 2mn) = 1 = (m^2 + n^2, m^2 - n^2)$; and by **Fact 4,** $(m^2 + n^2, (2mn)(m^2 - n^2)) = 1$

### Conclusion:

1) The heights $h_\beta$ and $h_\gamma$ are always integers.
2) The height $h_\alpha$ is always a rational number; it is an integer only in those nonprimitive pythagorean triangles with $\delta$ being a positive integer multiple of $m^2 + n^2$.

(iii) $\delta_\alpha = \left(\dfrac{2\delta(2mn)(m^2-n^2)}{2mn+m^2-n^2}\right)\cdot\dfrac{1}{\sqrt{2}}$

$\delta_\beta = \left(\dfrac{\delta\cdot(m^2+n^2)(m^2-n^2)}{m^2}\right)\cdot\sqrt{\dfrac{m}{m^2+n^2}}$,

$\delta_\gamma = \left(\dfrac{2\delta(2mn)(m^2+n^2)}{(m+n)^2}\right)\cdot\sqrt{\dfrac{m^2-n^2}{2(m^2+n^2)}}$

Looking at $\delta_\alpha$, we see that it is the product of a rational number with the irrational number $\dfrac{1}{\sqrt{2}}$; hence, $\delta_\alpha$ must be irrational.

Next consider $\delta_\gamma$. In a first course in elementary theory, a student is asked to show that the $n$th root of a natural number is rational if, and only if, that natural number is the $n$th power of another natural number; in particular the square root of a natural number is rational if, and only if, that natural number is a perfect square (which means that the square root will be a positive integer). This result may be found in number theory books. Two excellent sources are W. Sierpinski's book (see [1]) and Kenneth H. Rosen's *Elementary number Theory and Its Applications* (for details, refer to [6]). Taking this result a step further if $\gamma$ is a rational number $r=\dfrac{a}{b}$, where $a$ and $b$ are positive integers which are relatively prime; then $\sqrt{r}$ is rational if, and only if, both $a$ and $b$ are perfect squares. Applying this idea to $\delta_\gamma$ we see that $\delta_\gamma$ will be rational if, and only if, both numbers $(m^2-n^2)$ and $2(m^2+n^2)$ are integer squares; there is a detail that needs to be filled in here; namely that $(m^2-n^2, 2(m^2+n^2))=1$. This follows from $(m,n)=1, m+n\equiv 1(\mod 2)$, **Facts 4 and 6.**

However, $2(m^2+n^2)$ cannot be an integer square; if it were, then $2(m^2+n^2)=t^2$, for some positive integer; 2 must divide $t^2$; and thus (since it is a prime) 2 must divide t; set $t=2t_1$ to obtain $2(m^2+n^2)=4t_1^2$, or $m^2+n^2=2t_1^2$ which is a contradiction since $m^2+n^2$ is an odd number; while $2t_1^2$ is even.

This last argument on equation $2(m^2+n^2)=t^2$ could also be formulated in the language of congruences as follows: since $m^2+n^2\equiv 1(\mod 2)$, it follows $2(m^2+n^2)\equiv 2(\mod 4)$; but $t^2\equiv 0$ or 1(mod4) (depending on whether t is even or odd respectively), hence a contradiction.

In any event, since $2(m^2+n^2)$ can not be a perfect square the real number $\sqrt{\dfrac{m^2-n^2}{2(m^2+n^2)}}$ must be irrational for reasons already explained above. Hence $\delta_r$ must be irrational since it is the product of a rational number with an irrational number.

Now, let us look at $\delta_\beta$; since $(m, m^2+n^2)=1$ (even though not included among the seven listed facts, it is clear that if $(a,b)=1$, then $(a^{\varepsilon 1}, a^{\varepsilon 1}+b^{\varepsilon 2})=1$, for any positive integer exponents $\varepsilon_1$ and $\varepsilon_2$), repeating the first part of the analysis given for the irrationality of $\delta_\gamma$, we see that $\delta_\beta$ will be rational if, and only if, both numbers m and $m^2+n^2$ are integer squares; $m=t^2$ and $m^2+n^2=z^2$; for natural numbers $t$ and $z$.

What this really says is that the numbers $m$ and $n$ would be the leg lengths of a primitive pythagorean triple $(m,n,z)$ and with the leg length m being a perfect square.

The conditions $m^2+n^2=z^2, (m,n)=1, m+n\equiv 1(\mod 2)$ imply either

$$\left.\begin{array}{l}(m=t_1^2-t_2^2, n=2t_1t_2, z=t_1^2+t_2^2)\\ or\\ (m=2t_1t_2, n=t_1^2-t_2^2, z=t_1^2+t_2^2)\end{array}\right\} \qquad (2)$$

for some natural numbers $t_1, t_2$ such that $(t_1, t_2)=1, t_1>t_2$ and $t_1+t_2\equiv 1(\mod 2)$. Assume that the first set of equations holds in (2). If we impose the condition $m=t^2$, we see that, since $m$ is odd, we must have $t^2=t_1^2-t_2^2$. We end up with another primitive Pythagorean triple. $(t,t_2,t_1): t^2+t_2^2=t_1^2$. Resolving this we obtain,

$$(t=t_3^2-t_4^2, t_2=2t_3t_4, t_1=t_3^2+t_4^2)$$

for natural numbers $t_3, t_4$ with $(t_3, t_4)=1, t_3>t_4$, and $t_3+t_4\equiv 1(\mod 2)$.

Finally, if the second set of equations holds in (2), and we impose the condition $m = t^2$, some number theory shows that either $t_1 = 2t_5^2$ and $t_2 = t_6^2$ or vice-versa, namely $t_1 = t_5^2$ and $t_2 = 2t_6^2$: for natural numbers $t_5, t_6$ such that $(t_5, t_6) = 1$; and in first case with $t_6$ odd; while in the second case with $t_5$ odd. In both cases, $t = 2t_5 t_6$.

All in all, the resolution of the simultaneous conditions $m = t^2$ and $m^2 + n^2 = z^2$ leads to three families of types of numbers $m$ and $n$; below, in each case we use the letters $k$ and $l$ to describe these three two-parameter families. Note that if we substitute for $m = t^2$ and $m^2 + n^2 = z^2$ in the formula for $\delta_\beta$ we find,

$$\delta_\beta = \frac{\delta \cdot z \cdot (m^2 - n^2)}{t^3}; \text{ and by virtue of } (m, n) = 1,$$

$z^2 = m^2 + n^2$ and $_{m\ =\ t^2}$, one can easily show that $(z.(m^2 - n^2), t^3) = 1$; what this shows is that $\delta_\beta$ will be an integer if, $\delta$ is divisible by $t^3$. We have the following summary:

We mention to the reader that to express the first family listed below ( in terms of the independent integer parameters $k$ and $l$), one starts with $t_3 = k$ and $t_4 = l$; and then traces back to $t, t_1, t_2, m, n$ and z. In the case of families 2 and 3, one starts with $t_5 = k$ and $t_6 = l$; and then traces back to $t, t_1, t_2, m, n$ and z.

**Conclusion:** (By (1) ) Keep in mind that in all three families we must have $m > n$

1) The internal bisectors $\delta_\alpha$ and $\delta_\gamma$ are always irrational numbers.
2) The internal bisector of $\delta_\beta$ is irrational unless $\delta_\beta = \dfrac{\delta \cdot z(m^2 - n^2)}{t^3}$, where $\delta$ is a natural number; and the positive integers $z, m, n, t$ belong to any of the families described below; for each of the three families, the resulting value of $\delta_\beta$ is rational. The value of $\delta_\beta$ is integral if, and only if, $\delta$ is an integer multiple of $t^3$.

**Family 1:**
$$\begin{cases} m = (k^2 + l^2)^2 - (2kl)^2 = (k^2 - l^2)^2 \\ n = 2(k^2 + l^2)(2kl) = 4kl(k^2 + l^2) \\ z = (k^2 + l^2) + (2kl)^2 = k^4 + l^4 + 6k^2l^2 \\ t = k^2 - l^2 \end{cases}$$

Where the positive integers $k, l$ satisfy $(k, l) = 1, k > l,$ and $k + l \equiv 1 \pmod{2}$; and such that $m > n$.

**Family 2:**
$$\begin{cases} m = 4k^2l^2 \\ n = 4k^4 - l^4 \\ z = 4k^4 + l^4 \\ t = 2kl \end{cases}$$

Where the positive integers $k, l$ satisfy $(k, l) = 1, l$ odd, and $2k^2 > l^2$ (since n is positive); and such that $m > n$.

**Family 3:**
$$\begin{cases} m = 4k^2l^2 \\ n = k^4 - 4l^4 \\ z = k^4 + 4l^4 \\ t = 2kl \end{cases}$$

Where the positive integers $k, l$ satisfy $(k, l) = 1, k$ odd, and $k^2 > 2l^2$; and such that $m > n$

(iv)

$$d_\alpha = \left( \frac{2\delta(2mn)(m^2 - n^2)}{|2mn - (m^2 - n^2)|} \right) \cdot \frac{1}{\sqrt{2}}$$

$$d_\beta = \left( \frac{\delta(m^2 + n^2)(m^2 - n^2)}{n} \right) \cdot \frac{1}{\sqrt{m^2 + n^2}}$$

$$d_\gamma = \left( \frac{2\delta(2mn)(m^2 + n^2)}{m - n} \right) \cdot \frac{1}{\sqrt{2(m^2 + n^2)}}$$

Because of all the work, techniques, and arguments we used in the previous case, case (iii) (internal bisectors), our conclusions here are swift indeed: by inspection $d_\alpha$ is irrational; also, as we saw in the case of the internal bisectors, $2(m^2 + n^2)$ can not be an integer square, which establishes the irrationality of $d_\gamma$; finally, $d_\beta$ will be rational precisely when $m^2 + n^2 = w^2$, where $w$ is some natural number. This shows that $d_\beta$ will be rational exactly when $m$ and $n$ are the leg lengths of a primitive Pythagorean triple $(m, n, w)$. Which means either

$(m = K^2 - L^2, n = 2KL, w = K^2 + L^2)$

or $(m = 2KL), n = K^2 - L^2, w = K^2 + L^2)$,

for some positive integers $K, L$ such that $(K, L) = 1, K > L$, and $K + L \equiv 1 \pmod{2}$.

**Conclusion:**

1) Both external bisectors $d_\alpha$ and $d_\gamma$ are irrational numbers.
2) The external bisectors $d_\beta$ is irrational unless $d_\beta = \dfrac{\delta \cdot w \cdot (m^2 - n^2)}{n}$, where $\delta$ is a natural number; and the positive integers $w, m, n$ belong to either of the two families of numbers described below; for each of the two families, the resulting value of $d_\beta$ is rational. The value of $d_\beta$ is integral if, and only if, $\delta$ is an integer multiple of n (note that $(note\ that\ (n, w(m^2 - n^2)) = 1)$

**Family 4:** $\begin{cases} m = K^2 - L^2 \\ n = 2KL \\ w = K^2 + L^2 \end{cases}$, where the positive integers $K, L$

Satisfy $(K, L) = 1, K > L$ and $K + L \equiv 1 (\mod 2)$; and with $m > n$.

**Family 5:** $\begin{cases} m = 2KL \\ n = K^2 - L^2 \\ w = K^2 + L^2 \end{cases}$, where the positive integers $K, L$

Satisfy $(K, L) = 1, K > L$ and $K + L \equiv 1 (\mod 2)$; and with $m > n$.

(v) $\rho_\alpha = \delta n(m+n), \rho_\beta = \delta n(m+n), \rho_\gamma = \delta n(m-n)$. Obviously, all three radii are integers.

(vi) $\mu_\alpha = \dfrac{\delta(m^2 + n^2)}{2}$; $\mu_\alpha$ is half of an odd integer if $\delta$ is odd; while if $\delta$ is even, it is an integer.

# 5 The Irrationality of $\mu_\beta$ and $\mu_\gamma$

$$\mu_\beta = \delta.\sqrt{m^4 + n^4 - m^2 n^2}, \mu_\gamma = \dfrac{\delta.\sqrt{m^4 + n^4 + 14m^2 n^2}}{2}$$

We know that the natural numbers $m$ and $n$ satisfy the conditions $(m,n)=1$ and $m+n \equiv 1 \pmod{2}$. According to results (A) and (B) listed below, neither of the numbers $m^4+n^4-m^2n^2$ and $m^4+n^4+14m^2n^2$ can be a perfect square; which implies, in accordance with the two formulas above, that both medians $\mu_\beta$ and $\mu_\gamma$ must be irrational numbers.

(A) Under the conditions $(x,y)=1$ and $x+y \equiv 1 \pmod{2}$ (i.e. x odd, y even; or vice-versa), the three-variable diophantine equation $x^4-x^2y^2+y^4=z^2$ has no solution in positive integers; if we keep the condition $(x,y)=1$ but we drop the condition $x+y \equiv 1 \pmod{2}$ (so that both $x, y$ are odd), the same equation has a unique solution in positive integers, namely $x=y=z=1$. Finally, if we drop both of these conditions the general solution to the same equation in positive integers is given by $x=y=t, z=t^2$; when $t$ can be any positive integer.

(B) Under the conditions $(x,y)=1$ and $x+y \equiv 1 \pmod{2}$, the three variable Diophantine equation, $x^4+14x^2y^2+y^4=z^2$, has no solution in the set of positive integers (note, however, that if we replace the condition $x+y \equiv 1 \pmod{2}$ by the condition $x \equiv y \equiv 1 \pmod{2}$, there is at least one solution; $x=y==1, z=4$).

The first historically known proof of (A) is due to Pocklington (see [4]). The same proof can also be found in W. Sierpinski's number theory book (see [1]). Pocklington's proof is clever one, but certainly not an easy one. Below we offer a brief outline of Pocklinton's proof for the interested reader. Also, we present a proof, due to this author, of result (B). We point out here, that, in general, three-variable Diophantine equations of the form $ax^4+bx^2y^2+cy^4=z^2$, with the coefficients $a,b,c$ being and integers, are extremely difficult to deal with.

As a historical note, Euler is cited in Dickson's book (see [5] for details) as having proved, in 1780, result (B); also as having shown that the same Diophantine equation, under the conditions $(x,y)=1$ and with both $x$ and $y$ begin odd, has only one solution in positive integers; namely $(x,y,z)=(1,1,4)$.

(Please note: This author has neither seen nor is familiar with Euler's proof).

## An outline of Pocklington's proof of result (A)

The proof starts by pointing out that if the equation $x^4 - x^2y^2 + y^4 = z^2$ is solvable in positive integers $x, y, z$ with $(x, y) = 1$ and $x + y \equiv 1 \pmod{2}$, then by Fact 1 there must exist a positive integer solution $(x_0, y_0, z_0)$ to equation (A) such that $x_0 + y_0 \equiv 1 \pmod{2}, (x_0, y_0) = 1$, and with the product $x_0 y_0$ being least; then $x_0^4 - x_0^2 y_0^2 + y_0^4 = z_0^2$ implies $(x_0^2 - y_0^2)^2 + (x_0 y_0)^2 = z_0^2$; thus, since the condition $(x_0, y_0) = 1$ easily implies $(x_0^2 - y_0^2, x_0 y_0) = 1$, it follows that the triple $(x_0^2 - y_0^2, x_0 y_0, z_0)$ is a primitive one with the integer $x_0^2 - y_0^2$ being odd and $x_0 y_0$ even. Consequently, we must have $x_0^2 - y_0^2 = M^2 - N^2, x_0 y_0 = 2MN, z_0 = M^2 + N^2$, for some positive integers $M, N$ with $(M, N) = 1, M + N \equiv 1 \pmod{2}$ and $M > N$ (because of the symmetry of the initial equation with respect to the variables $x, y$; and $x \neq y$ (which follows from the condition $x + y \equiv 1 \pmod{2}$ ), we may assume $x_0 > y_0$. At this point, the proof branches out a bit; in the end by using some basic results from number theory including Facts 2 and 5; a new solution is $(x_1, y_1, z_1)$ is found, which satisfies $x_1 y_1 < x_0 y_0$; this contradicts that $x_0 y_0$ is minimal.

The above description of proof is an example of the method of *infinite descent*, first established by P. Fermat. In a nutshell, this method goes like this: If $f(x_1, x_2, ..., x_n) = 0$ is a Diophantine equation (here $f$ stands for some polynomial function in $n$ variable; with integer coefficients), then by assuming the existence of a solution $(a_1, a_2, ..., a_n)$; in positive integers $a_1, a_2, ..., a_n$; one shows that this leads to another solution $(b_1, b_2, ..., b_n)$; also in positive integers $b_1, b_2, ..., b_n$. Also, one establishes that $g(b_1, b_2, ..., b_n) < g(a_1, a_2, ..., a_n)$, where $g$ is a well-defined algebraic function (in $n$ variables), which, when restricted over the set of positive integers, its range is a subset of $Z^+$ (set of positive integers). This clearly violates Fact 1 (well ordering principle), for one can then find yet another solution $(c_1, c_2, ..., c_n)$ such that $g(c_1, c_2, ..., c_n) < g(b_1, b_2, ..., b_n)$ and so on, *ad infinitum*. At some point in our proof of result (B): we employ the infinite descent method.

*Proof of result* (B): Let $(a, b, c)$ be a solution in positive integers of the equation in (B); we have,

$$(a^2 - b^2)^2 + (4ab)^2 = c^2$$
$$(a, b) = 1, a + b \equiv 1 \pmod{2}$$
(3)

Because the Diophantine equation in (B) has symmetry in the variables $x, y$ and since $x \neq y$ (by virtue of $x + y \equiv 1 \pmod{2}$), we may assume $a > b$. The conditions $(a,b) = 1$ and $a + b \equiv 1 \pmod{2}$, imply $(a^2 - b^2, 4ab) = 1$; thus the triple $(a^2 - b^2, 4ab, c)$ is a primitive Pythagorean one with $a^2 - b^2$ being odd, while $4ab$ is even. We must have

$$a^2 - b^2 = M^2 - N^2 \qquad (4_1)$$

$$\text{and } 4ab = 2MN \qquad (4_2)$$

with the positive integers $M, N$ satisfying $(M, N) = 1, M + N \equiv 1 \pmod{2}$, and $M > N$. Let us write equation $(4_2)$, in the form,

$$\frac{2a}{M} = \frac{N}{b} \qquad (5)$$

We set $\dfrac{a}{M} = r$, a positive rational number. From (5) we obtain $a = M \cdot r$ and $N = 2rb$; combining these two with equation $(4_1)$ and with a bit of algebra we arrive at,

$$\left(\frac{M}{b}\right)^2 = \frac{4r^2 - 1}{1 - r^2} \qquad (6)$$

Obviously, according to (6), either both rational numbers $4r^2 - 1$ and $1 - r^2$ are positive; or they are both negative; but $0 < r$, so if it were $4r^2 - 1 < 0$ and $1 - r^2 < 0$; we would then have $0 < r < \dfrac{1}{2}$ and simultaneously, $1 < r$; an impossibility. Therefore, we must have $4r^2 - 1 > 0$ and $1 - r^2 > 0$ (which imply, by virtue of $r > 0$, $\dfrac{1}{2} < r < 1$). Now we write the rational number $r$ in lowest terms: $r = \dfrac{t_1}{t_2}$, with the positive integers $t_1, t_2$ being relatively prime: $(t_1, t_2) = 1$. We substitute for $r = \dfrac{t_1}{t_2}$ in equation (6) to obtain,

$$\left(\frac{M}{b}\right)^2 = \frac{4t_1^2 - t_2^2}{t_2^2 - t_1^2} \tag{7}$$

Also mote that since $\frac{1}{2} < r < 1$ (see above) we must have $1 \leq t_1 < t_2$ and $2 \leq t_2 < 2t_1$. Now let $d$ be that greatest common divisor of the positive integers $4t_1^2 - t_2^2$ and $t_2^2 - t_1^2$; $d = \left(4t_1^2 - t_2^2, t_2^2 - t_1^2\right)$. According to the **Lemma** (which can be found at the end of this proof), since the rational number $\frac{4t_1^2 - t_2^2}{t_2^2 - t_1^2}$ is a rational square (equation (7) ), each of the numbers $4t_1^2 - t_2^2$ and $t_2^2 - t_1^2$ must be equal to $d$ times an integer square:

$$\left. \begin{array}{l} 4t_1^2 - t_2^2 = d.u^2 \\ t_2^2 - t_1^2 = d.v^2 \end{array} \right\} \tag{8}$$

Here the positive integers $u$ and $v$ are relatively prime: $(u,v) = 1$. What are the possible values of $d$? The answer is that $d$ can only be equal to 1 or 3. This is easy to see: by adding the two equations in (8) we see that $d$ must be a divisor of $3t_1^2$; and by taking the sum of the first equation in (8) plus four times the second equation, we see that $d$ must also divide $3t_2^2$. In effect, $d$ must be a common divisor of $3t_1^2$ and $3t_2^2$; but $(t_1, t_2) = 1$, which easily implies $d = 1$ or 3. The possibility $d = 1$ is easily eliminated: indeed, if we see $d = 1$ in (8) we obtain,

$$4t_1^2 - t_2^2 = u^2$$

$$t_2^2 - t_1^2 = v^2; t_2^2 = t_1^2 + v^2$$

Since $(t_2, t_1) = 1$ ( and consequently $(t_2, v) = 1 = (t_1, v)$), we conclude form the second equation that ( $t_1, v, t_2$) is a primitive pythagorean triple, so that $t_2$ must be odd; which shows $t_2^2 \equiv 1 \pmod 4$. On the other hand, the first equation above implies that $u$ must be odd, since $t_2$ is odd and $4t_1^2$ is even. Hence $t_2 \equiv u \equiv 1 \pmod 2 \Rightarrow t_2^2 \equiv u^2 \equiv 1 \pmod 4$. But, $u^2 + t_2^2 = 4t_1^2 \equiv 0 \pmod 4$, while $u^2 + t_2^2 \equiv 1 + 1 \equiv 2 \pmod 4$, a contradiction.

Now consider the only other possibility remaining: $d = 3$. From equations (8) we obtain

$$\left.\begin{array}{l}4t_1^2 - t_2^2 = 3u^2 \\ \\ t_2^2 - t_1^2 = 3v^2\end{array}\right\}; \text{ solving for } t_1^2 \text{ and } t_2^2 \text{ gives,}$$

$$\left.\begin{array}{l}t_1^2 = u^2 + v^2 \\ \\ t_2^2 = u^2 + (2v)^2\end{array}\right\} \tag{9}$$

First note that $u$ and $v$ can not both be odd: for in that case we would have $u^2 \equiv v^2 \equiv 1 (\mod 4)$ and the first equation in (9) would imply $t_1^2 \equiv 2(\mod 4)$ which is impossible. Therefore, in view of $(u,v) = 1$, $u$ and $v$ must have different parity.

There are two possibilities remaining: either $u$ is even and $v$ is odd, or alternatively, $u$ is odd and $v$ is even.

Let us first assume that $u$ is even and $v$ is odd. Set $u = 2u_1$, where $u_1$ is a positive integer. From equations (9) we obtain $t_2 = 2t_3$ and,

$$\left.\begin{array}{l}t_1^2 = (2u_1)^2 + v^2 \\ \\ t_3^2 = u_1^2 + v^2\end{array}\right\} \tag{10}$$

Bu virtue of $(u,v) = 1, u = 2u_1$, and $v$ being odd, we easily see that $(u_1, v) = 1 = (2u_1, v)$. Thus, according to (10) both $(2u_1, v, t_1)$ and $(u_1, v, t_3)$ are primitive pythgorean triples and so we must have (since $v$ is odd),

$$\left.\begin{array}{l}v = M_1^2 - N_1^2, 2u_1 = 2M_1 N_1, t_1 = M_1^2 + N_1^2 \\ \text{and} \\ v = M_3^2 - N_3^2, u_1 = 2M_3 N_3, t_3 = M_3^2 + N_3^2\end{array}\right\} \tag{11}$$

for positive integers $M_1, N_1, M_3, N_3$ satisfying $(M_1, N_1) = 1 = (M_3, N_3)$ and $M_1 + N_1 \equiv 1 \equiv M_3 + N_3 \pmod{2}$.

Using (the identity) $(M_1^2 - N_1^2)^2 + (2M_1 N_1)^2 = (M_1^2 + N_1^2)^2$ and equations (11) and (10) we arrive at,

$$(M_3^2 - N_3^2)^2 + (4M_3 N_3)^2 = t_1^2;$$

$$M_3^4 + 14 M_3^2 N_3^2 + N_3^4 = t_1^2 \tag{12}$$

According to (12) the triple $(M_3, N_3, t_1)$ is a positive integer solution to the Diophantine equation in (B); by (7) and (9) we have,

$$\left(\frac{M}{b}\right)^2 = \frac{4t_1^2 - t_2^2}{t_2^2 - t_1^2} = \frac{3u^2}{3v^2} = \frac{u^2}{v^2}$$

The relation $\dfrac{M^2}{b^2} = \dfrac{u^2}{v^2}$, in conjunction with $(u,v) = 1$, shows that $M^2 \geq u^2$ and $b^2 \geq v^2$ (consider the greater common divisor of $M$ and $b$). Thus, $M^2 \geq u^2 = (2u_1)^2 > u_1^2$, and since $N^2 M^2 \geq M^2$, we conclude that

$N^2 M^2 > u_1^2 \Rightarrow$ (since $NM$ and $u_1$ are positive) $NM > u_1 \Rightarrow$ (by (11)) $2NM > 4M_3 N_3$, and by ($4_2$) we arrive at $4ab > 4M_3 N_3; ab > M_3 N_3$. So both $(a,b,c)$ and $(M_3, N_3, t_1)$ are positive integer solutions to the Diophantine equation in (B), and they satisfy $M_3 N_3 < ab$.

Below we show that the same phenomenon occurs if in (9) we assume $u$ to be odd and $v$ to be even. In view of $(v, u) = 1$ we must have $(2v, u) = 1$; therefore (9) shows that both $(u, v, t_1)$ and $(u, 2v, t_2)$ are primitive pythagorean triples:

$$\left.\begin{array}{l}u = m_1^2 - n_1^2,\ v = 2m_1 n_1,\ t_1 = m_1^2 + n_1^2 \\ u = m_2^2 - n_2^2,\ 2v = 2m_2 n_2,\ t_2 = m_2^2 + n_2^2\end{array}\right\} \quad (13)$$

for positive integers $m_1, n_1, m_2, n_2$ satisfying $m_1 + n_1 \equiv 1 \equiv m_2 + n_2 \pmod{2}$, and $(m_1, n_1) = 1 = (m_2, n_2)$. Applying the identity, $(m_2^2 - n_2^2)^2 + (2m_2 n_2)^2 = (m_2^2 + n_2^2)^2$; and using (13) we arrive at $(m_1^2 - n_1^2)^2 + (4m_1 n_1)^2 = t_2^2$;

$$m_1^4 + 14 m_1^2 n_1^2 + n_1^4 = t_2^2. \quad (14)$$

Equation (14) clearly shows that $(m_1, n_1, t_2)$ is a positive integer solution to the diophantine equation in (B). Again, as before, using $\dfrac{M^2}{b^2} = \dfrac{u^2}{v^2}$ (it follows from (7) and (9)) and $(u,v) = 1$; we infer $M^2 \geq u^2$ and $b^2 \geq v^2$.

Clearly, $ab \geq b \geq v = 2m_1 n_1 > m_1 n_1$. Thus both $(a,b,c)$ and $(m_1, n_1, t_2)$ are positive integer solutions to the diophantine equation in (B), and they satisfy $m_1 n_1 < ab$. In summary, we see that in both cases (i.e., $u$ even, $v$ odd; and $u$ odd, $v$ even), the two equations in (9) lead to a new positive integer solution $(a', b', c')$ to the Diophantine equation in (B), such that $a'b' < ab$. This, though, violates **Fact 1**, which, in effect implies that the set of all positive integer solutions to the equation in (B), must have a solution $(a_0, b_0, c_0)$ with $a_0 b_0$ being least. (The process could be repeated *ad infinitum,* each time producing a "smaller" solution in the sense described above.)

The proof is now complete: the assumption that the Diophantine equation in (B) has a positive integer solution leads to a contradiction; thus, it has no such solutions.

**Lemma:** If $n$ is a natural number, and the positive integers $a,b,c,d$ satisfy $\dfrac{a}{b} = \dfrac{c^n}{d^n}$, then $a = \delta \cdot c_1^n$ and $b = \delta \cdot d_1^n$, where $\delta = (a,b)$ and $a_1, b_1$ are relatively prime positive integers; $(a_1, b_1) = 1$

**Proof:** We put $a = \delta \cdot \alpha$, $b = \delta \cdot \beta$, where $\alpha, \beta$ are natural numbers such that $(\alpha, \beta) = 1$

We have,

$$\frac{\delta \cdot \alpha}{\delta \cdot \beta} = \frac{c^n}{d^n} \Rightarrow \alpha \cdot d^n = \beta \cdot c^n$$

Let $D = (d, c)$, so that $d = D \cdot d_1, c = D \cdot c_1$, for positive integers $d_1, c_1$ with $(c_1, d_1) = 1$. Thus, $\alpha \cdot D^n \cdot d_1^n = \beta \cdot D^n \cdot c_1^n \Rightarrow \alpha \cdot d_1^n = \beta \cdot c_1^n$.

According to **Fact 4,** the relation $(c_1, d_1) = 1$ implies $(c_1^n, d_1^n) = 1$. Now look at $\alpha \cdot d_1^n = \beta \cdot c_1^n$. Since $(c_1^n, d_1^n) = 1 = (\alpha, \beta)$, Fact 3 implies that $\alpha$ is a divisor of $c_1^n$ and $c_1^n$ is a divisor of $\alpha$. Hence by Fact 5, we must have $\alpha = c_1^n$; and therefore $d_1^n = \beta$ as well. Therefore, $a = \delta \cdot \alpha = \delta \cdot c_1^n$; and $b = \delta \cdot \beta = \delta \cdot d_1^n$. End of proof.

## 6 Numerical Examples

**Family 1:** The smallest value of $m$ is $m = (6^2 - 1^2) = (35)^2 = 1225$, obtained for $k = 6$ and $l = 1$. Correspondingly,

$$n = 4 \cdot 6 \cdot 1 \cdot (6^2 - 1^2) = 888,$$

$$z = 6^4 + 1^4 + 6 \cdot 6^2 \cdot 1^2 = 1513.$$

We have that $m^2 - n^2 = (1225)^2 - (888)^2 = (337)(2113) = 712,081$ and $t = 6^2 - 1^2 = 35$. Thus, $\delta_\beta = \frac{\delta \cdot (1513)(712,081)}{(35)^3} = \frac{1077378553\delta}{42875}$, and since 1077378553 is divisible by neither 5 nor by 7, we see that the fraction $\frac{1077378553}{42875}$ is reduced to lowest terms. For $\delta = 1$, the primitive pythagorean triangle with sides $\beta = 2mn = 2175600, \gamma = m^2 - n^2 = 712,081, \alpha = m^2 + n^2 = 2289169$, has internal bisector $\delta_\beta = \frac{1077378553}{42875}$

**Family 2:** the smallest value of $m$ is 1600 it is obtained for $k = 4$ and
$l = 5 : m = 4 \cdot 4^2 \cdot 5^2 = 1600, n = 4 \cdot 4^4 - 5^4 = 1024 - 625 = 399, z = 4 \cdot 4^5 + 5^4 = 1024 + 625 = 1649,$
$t = 2 \cdot 4 \cdot 5 = 40$

We have $\delta_\beta = \dfrac{\delta \cdot 1649 \cdot (1600 - 399)(1600 + 399)}{(40)^3} = \dfrac{\delta \cdot 1649 \cdot 2400799}{(40)^3} \Rightarrow \delta_\beta = \dfrac{3958917551\delta}{(40)^3}$ so

for $\delta = 1$, we obtain the primitive pythagorean triangle with sides
$\beta = 2mn = 1276800, \gamma = m^2 - n^2 = 2400799$, and with $\delta_\beta = \dfrac{3958917551}{6400}$.

**Family 3:** Smallest value of $m$ is 144, obtained for $k = 3$ and $l = 2$. We have
$m = 4 \cdot 3^2 \cdot 2^2 = 144, n = 3^4 = 4 \cdot 2^4 = 17, z = 3^4.4.2^4 = 145, t = 2 \cdot 3 \cdot 2 = 12$; and thus,

$\delta_\beta = \dfrac{\delta \cdot 145 \cdot (144 - 17)(144 + 17)}{(12)^3} = \dfrac{\delta \cdot 145 \cdot 127 \cdot 161}{(12)^3} = \dfrac{2964815}{1728}$

For $\delta = 1$, we obtain the primitive pythagorean triangle with sides
$\beta = 2mn = 4896 = \gamma = m^2 - n^2 = 20447, a = 21025$, and with $\delta_\beta = \dfrac{2964815}{1728}$.

**Family 4:** For $K = 4$ and $L = 1$, we have $m = 4^2 - 1^2 = 15, n = 2 \cdot 4 \cdot 1 = 8, w = 4^2 + 1^2 = 17$,
$d_\beta = \dfrac{\delta.17.(15 - 8)(15 + 8)}{8} = \dfrac{17.7.23\delta}{8} = \dfrac{2737\delta}{8}$.

Now, the smallest value of $m$ is 3; obtained for $K = 2$ and $L = 1$; which gives
$m = 3, n = 4, w = 5, B = 24\delta, \gamma = 7\delta, \alpha = 25\delta, d_\beta = \dfrac{35\delta}{4}$. Note that for $\delta = 4$ we obtain Diophantus' example (see introduction), namely $\beta = 96, \gamma = 28, \alpha = 100, d_\beta = 35$. For $\delta = 1$, we obtain the primitive pythagorean triangle with side lengths
$\beta = 2mn = 2 \cdot 15 \cdot 8 = 240, \gamma = m^2 - n^2 = 161, \alpha = m^2 + n^2 = 289$, and external bisector $d_\beta = \dfrac{2737}{8}$.

**Family 5:** the smallest value of $m$ is 4 obtained from $K = 2$ and $L = 1$; we have
$m = 2 \cdot 2 \cdot 1 = 4, n = 2^2 - 1^2 = 3, w = 2^2 + 1^2 = 5, d_\beta = \dfrac{\delta.5(4 - 3)(4 + 3)}{3} = \dfrac{35\delta}{3}$. For $\delta = 1$, we obtain the primitive pythaforean triangle with side lengths
$\beta = 2mn = 2 \cdot 4 \cdot 3 = 24, \gamma = m^2 - n^2 = 7, \alpha = m^2 + n^2 = 25$, and external bisector $\delta_\beta = \dfrac{35}{3}$.

The derivation and presentation of the various formulas for the seventeen secondary elements of pythagorean triangles can be found on pages.

The formulas for the internal angle bisectors and the medians can be found in Chapter 5, starting with Section 5.1.